\documentclass{ntmanuscript}
\title{Facility Deployment Decisions through Warp Optimizaton 
       of Regressed Gaussian Processes}

\author{Anthony Michael Scopatz$^1$}

\institute{$^1$University of South Carolina, Department of Mechanical 
    Engineering, Nuclear Engineering Program, Columbia, SC 29201}

\submitter{Anthony M. Scopatz}
\submitteraddress{541 Main Street, Columbia, SC 29208}
\submitteremail{scopatz@cec.sc.edu}


\usepackage{color}
\usepackage{graphicx}
\usepackage{booktabs} 
\usepackage{microtype} 
\usepackage{xspace}
\usepackage{listings}
\usepackage{textcomp}
\usepackage{ulem}
\usepackage{amssymb}
\DeclareMathAlphabet{\mathpzc}{OT1}{pzc}{m}{it}

\definecolor{listinggray}{gray}{0.9}
\definecolor{lbcolor}{rgb}{0.9,0.9,0.9}
\newcommand{\E}{\mathbb{E}}
\newcommand{\GP}{\mathpzc{GP}}
\newcommand{\N}{\mathbb{N}}

\newcommand{\argmin}{\mathrm{argmin}}

\newcommand{\I}{\mathbf{I}}
\newcommand{\K}{\mathbf{K}}
\newcommand{\stochastic}{\texttt{`stochastic'}\xspace}
\newcommand{\innerprod}{\texttt{`inner-prod'}\xspace}
\newcommand{\allflag}{\texttt{`all'}\xspace}

\date{}
\begin{document}

\begin{abstract}
A method for quickly determining deployment schedules that meet a given 
fuel cycle demand is presented here. This algorithm is fast enough to 
perform \emph{in situ} within low-fidelity fuel cycle simulators. It uses
Gaussian process regression models to predict the production curve as a 
function of time and the number of deployed facilities. Each of these
predictions is measured against the demand curve using the dynamic time
warping distance. The minimum distance deployment schedule is evaluated
in a full fuel cycle simulation, whose generated production curve 
then informs the model on the next optimization iteration. The method
converges within five to ten iterations to a distance that is less than one 
percent of the total deployable production. A representative once-through
fuel cycle is used to demonstrate the methodology for reactor deployment.

\end{abstract}

\section{Introduction}
\label{intro}

With the recent advent of agent-based nuclear fuel cycle simulators, such as 
Cyclus \cite{DBLP:journals/corr/HuffGCFMOSSW15,cyclus_v1_0}, there comes the 
possibility to make \emph{in situ}, dynamic facility deployment decisions.
This  would more fully model real-world fuel cycles where institutions 
(such as utility companies)
predict future demand and choose their future deployment schedules 
appropriately. However, one of the major challenges to making \emph{in situ}
deployment decisions is the speed at which ``good enough'' decisions can 
be made. This paper proposes three related deployment-specific optimization 
algorithms that can be used for any demand curve and facility type.

The demands of a fuel cycle scenario can often be simply stated, e.g. 
1\% growth in power production [GWe]. Picking a deployment schedule for a 
certain kind of facility (e.g. reactors) can thus be seen as an optimization 
problem of how well the deployment schedule meets the demand. Here, the 
dynamic time warping (DTW) \cite{muller} distance is minimized 
between the demand curve and the regression of a Gaussian Process model (GP) 
\cite{rasmussen2006gaussian} of prior simulations. This minimization produces
a guess for a deployment schedule which is subsequently tested using 
an actual simulator. This process is repeated until an optimal deployment
schedule for the given demand is found.

Importantly, by using the Gaussian process surrogates, the number of 
simulation realizations that must be executed as part of the optimization may 
be reduced to only a handful. Furthermore, it is at least two 
orders-of-magnitude faster to test the model than it is to run a single
low-fidelity fuel cycle simulation. Because of the relative computational 
cheapness, it 
is suitable to be used inside of a fuel cycle simulation. Traditional
\emph{ex situ} optimizers may be able to find more precise solutions but at a
computational cost beyond the scope and need of an \emph{in situ} use case
that is capable of dynamic adjustment.

Every iteration of the warp optimization of regressed Gaussian processes (WORG) 
method described here has two phases. The first is an estimation phase where 
the Gaussian process model is built and evaluated. The second takes the 
deployment schedule from the estimation phase and runs it in a fuel cycle 
simulator. The results of the simulator of the $s$-th iteration are then 
used to inform the model on the $(s+1)$-th iteration. 

Inside of each estimation phase there are three possible strategies for 
choosing the next deployment schedule.  The first is to sample of the 
space of all possible deployment strategies stochastically and then take the 
best guess.  The second is to search through the inner product of all choices,
picking the best option for each deployment parameter. The third strategy
is to perform the two previous strategies and determine which one has picked
the better guess.

Nuclear fuel cycle demand curve optimization faces many challenges. 
Foremost among these is that even though the demand curve is specified on 
the range of the real numbers, the optimization parameters are fundamentally 
integral in nature. For a discrete time simulator, deployments can only 
be issued in multiples of the size of the time step 
\cite{kelton2000simulation}. Furthermore, 
it is not possible to deploy only part of a facility; the facility is either 
deployed or it is not. While it may be possible to deploy a facility and 
only run it at partial capacity, most fuel cycle models do not support such
a feature for
keystone facilities.  For example, it is unlikely that a utility would build 
a new reactor only to run it at 50\% power. Thus, deployment is an integer 
programming problem, as opposed to its easier linear programming cousin
\cite{vanderbei2001linear}.

As an integer programming problem, the option space is combinatorially 
large. Assuming a 50 year deployment schedule where no more than 3 facilities 
are allowed to be deployed each time step, there are more than $10^{30}$ 
combinations. If every simulation took a very generous 1 sec, simulating 
each option would still take $\approx 3\times10^{12}$ times the current age
of the universe.

Moreover because all of the parameters are integral, there is not a 
meaningful formulation of a continuous Jacobian. Derivative-free optimizers are 
required. Methods such as particle swarm \cite{kennedy2010particle}, 
pswarm \cite{vaz2009pswarm}, and the 
simplex method \cite{vanderbei2001linear} all could work.  However, typical 
implementations require
more evaluations of the objective function (i.e. fuel cycle simulations)
than are within an \emph{in situ} budget. 

Even the usual case of 
Gaussian process optimization (sometimes known as kriging) 
\cite{osborne2009gaussian,simpson2001kriging} will still 
require too many full realizations in order to form an accurate model.
WORG, on the other hand, uses the dynamic time warping distance as a 
measure of how two time series differ. This is because the DTW distance is 
more separative than the typical
$L_1$ norm. Such additional separation drives the estimation phase to make 
better choices
sooner. This in turn helps the overall algorithm converge on a reasonable 
deployment schedule sooner. 
The stochastic strategy for WORG additionally utilizes Gaussian processes to 
weight the choice of parameters. This guides the guesses for the deployment
schedules such that fewer guesses are needed while simultaneously 
not forbidding  
any option.  So while WORG relies on Gaussian processes, it does so in a way
that is distinct from normal kriging. WORG
takes advantage of the \emph{a priori} knowledge that a deployment 
schedule is requested to meet a demand curve. This is not a strategy a 
generic, off-the-shelf optimizer would be capable of implementing.

The structure of the WORG algorithm is detailed in \S\ref{method}. 
The different strategies for selecting a best guess estimate of the 
deployment schedule are then discussed in \S\ref{selecting}. Performance
and results of the method for a sample once-through fuel cycle scenario 
are presented in \S\ref{results}. Finally, \S\ref{conclusion} summarizes
WORG and lists opportunities for future work.

\section{The WORG Method}
\label{method}

In order to describe the WORG method, first it is is useful to define  
notation for demand curves and their parameterization.  Call $t$ the time
[years] up to some maximal time horizon $T$ (e.g. 50 years) over which time 
the demand curve is known.  Then call $f(t)$ the demand curve in the natural 
units of the facility type (such as [GWe] for reactors).
$f(t)$ may be any function that is desired, including non-differential 
functions. For example, though, the demand curve for a 1\% growth rate 
starting at 90 [GWe] has the following form:
\begin{equation}
\label{f-1}
f(t) = 90\times 1.01^t
\end{equation}
Additionally, call $\Theta$ the deployment schedule for the facilities that 
may be constructed to meet the demand.
$\Theta$ is a sequence of $P$ parameters, indexed by $p$, as seen in 
Equation \ref{Theta}.
\begin{equation}
\label{Theta}
\Theta = \left\{\theta_1, \theta_2, \ldots, \theta_P\right\}
\end{equation}
Each $\theta_p$ represents that number of facilities to deploy on its
time step. In simple cases where there is only one type of facility
to deploy $P == T$.  However, when the deployment schedules of multiple 
facility types are needed to meet the same demand curve, $P > T$.  The usual
example for $P > T$ is for transition scenarios which necessarily require 
multiple kinds of reactors.

Now denote $M$ as the sequence for the minimum number of facilities deployable
for each deployment parameter. Also, call $N$ the sequence of the maximum number
of facilities deployable. The deployment parameters are thus each defined
on the range $\theta_p \in [M_p, N_p]$. Furthermore, because only whole
numbers of facilities may be deployed $\theta_p \in \N$.  It is also typical, 
but not required, for $M = \mathbf{0}$. Zero is also the lower bound
for all possible $\theta_p$ as facilities may not be forcibly retired via the
deployment schedule.

From here, call $g(t, \Theta)$ the production as a function of time for a
given deployment schedule. This has the same units as the demand curve.
Thus for power demand and reactor deployments, $g$ is in units of [GWe]. The 
optimization problem can now be posed as an attempt to find a $\Theta$
that minimizes the difference between $f$ and $g$.

\subsection{Dynamic Time Warping}
\label{dtw}

The question of how to take the difference between the demand curve and 
the production curve is an important one. The na\"ive option is to simply 
take the $L_1$ norm of the difference between these two time series, as 
seen in Equation \ref{delta-l1}.  However, since the $g(t, \Theta)$ computed
from a simulation is expensive, any operation that can meaningfully 
exacerbate the difference between time series helps drive down the number 
of optimization iterations.

Dynamic time warping is just such a mechanism. It computes 
a distance between any two time series which compounds the separation 
between the two. Additionally, the time series are not required to be of the 
same length, though for optimization purposes there is no reason for them 
not to be. DTW gives a measure of the amount that one time series would need to 
be warped to become the other time series. It is, therefore, a holistic  
measure that operates over all times. Dynamic time warping
is more fully covered in \cite{muller}.  However, an 
optimization-relevant introduction is given here.

For the time series $f$ and $g$, there are three parts to dynamic time 
warping. The first is the distance $d$, which will be minimized. The second 
is a cost matrix $C$ that helps compute $d$ by indicating how far a point 
on $f$ is from another point on $g$. Thirdly, the warp path $u$ is the 
minimal cost curve through the $C$ matrix from the fist point in time to 
the last. The DTW distance can thus be interpreted as the 
total cost of traveling the warp path.

The first step in computing a dynamic time warp distance is to 
assemble the cost matrix. Say that the demand time series $f$ has 
length $A$ indexed by $a$, and the production time series $g$ has 
length $B$ indexed by $b$. For the optimization problem here, $A$ and $B$
are in practice both equal to $T$.  However, it is useful to have $a$ and 
$b$ index the two time series separately. Now denote an $A\times B$ matrix 
$\Delta L$ as the $L_1$ norm of the difference between $f$ and $g$:
\begin{equation}
\label{delta-l1}
\Delta L_{a,b} = \left|f(a) - g(b, \Theta)\right|_1
\end{equation}
The cost matrix $C$ may now be defined as the $A\times B$ sized matrix 
which follows the recursion relations seen in Equation \ref{cost-matrix}.
\begin{equation}
\label{cost-matrix}
\begin{split}
C_{1,1} & = \Delta L_{1,1}\\
C_{1,b+1} & = \Delta L_{1,b} + C_{1,b}\\
C_{a+1,1} & = \Delta L_{a,1} + C_{a,1}\\
C_{a+1,b+1} & = \Delta L_{a,b} + \min\left[C_{a,b}, C_{a+1,b}, C_{a,b+1}\right]
\end{split}
\end{equation}
The boundary conditions above are the same as setting an infinite cost to 
any $a \le 0$ or $b \le 0$. The cost matrix $C$ has the same units as the 
demand curve. However, the scale of $C$ is 
larger than the demand, 
except for in the fiducial case. This is because the cost matrix compounds the 
minimum value of previous entries. 

Knowing a cost matrix, the warp path can be computed by traversing the 
matrix backwards from the $(A, B)$ corner to the $(1, 1)$ corner.
If the length of the warp is $I$ indexed by $i$, the warp path itself 
can be thought of as a sequence of coordinate points $u_i$. For a given 
point $u_i$ in the warp path, the previous point $u_{i-1}$ may be found by 
picking the minimum cost point among the locations one column over $(a,b-1)$, 
one row over $(a-1,b)$, and one previous diagonal element to $(a-1,b-1)$. 
Equation \ref{warp-path} expresses this mathematically.
\begin{equation}
\label{warp-path}
u_{i-1} = \argmin\left[C_{a-1,b-1}, C_{a-1,b}, C_{a,b-1}\right]
\end{equation}
The maximum possible length of $u$ is thus $\max(I) = A + B$.
The minimum possible length, though, is $\min(I) = \sqrt{A^2 + B^2}$. 

The dynamic time warping distance distance $d$ can now be stated as the 
cost of the final entry of the warp path normalized by the maximum possible
length of the warp path.  
\begin{equation}
\label{d-calc-ab}
d(f, g) = \frac{C_{A,B}}{A + B}
\end{equation}
However, because the demand curve and the production curve 
are often defined on the same time grid, $d$ can be further
reduced to the following:
\begin{equation}
\label{d-calc}
d(f, g) = \frac{C_{T,T}}{2T}
\end{equation}
Therefore, $d$ has the same units as the demand curve, production curve, 
and cost matrix.

\begin{figure}[htb]
\centering
\includegraphics[width=0.9\textwidth]{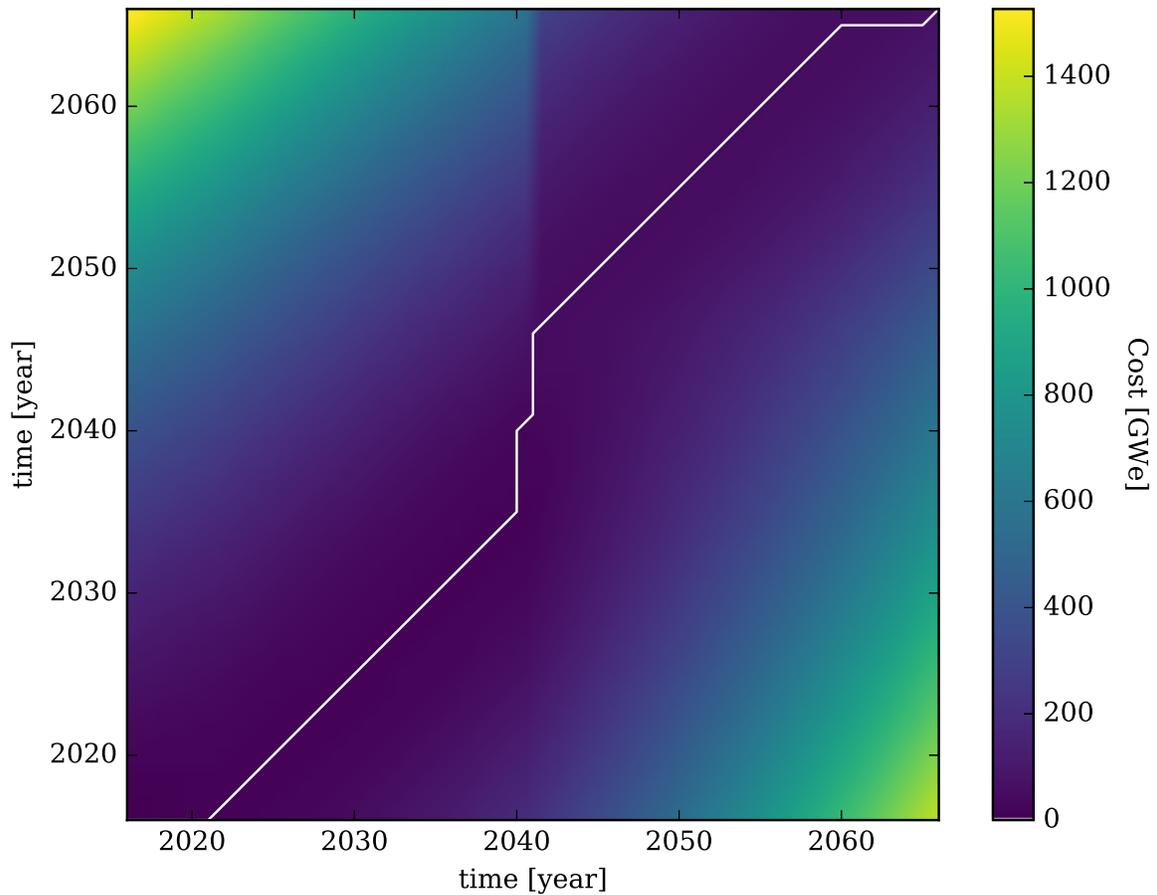}
\caption{Heat map of the cost matrix between a 1\% growth demand curve and 
a production curve the under produces by 5\% for the first 25 years and then
over produces for the second 25 years.  
The warp path $u$ is superimposed as the white curve on top of the 
cost matrix.}
\label{cost-demand-to-production}
\end{figure}

As an example, take a 1\% growth that starts with 90 GWe in the year 
2016 as the demand curve. Then consider a production curve that 
under-produces the demand by 5\% for 25 years before switching to 
over-producing this curve by 5\% for the next 25 years.  
Figure \ref{cost-demand-to-production} shows the dynamic time warping 
cost matrix between these two time series as a heat map.

Additionally, the warp path between the example demand and production 
curves is presented as the white curve on top of the heat map in 
Figure \ref{cost-demand-to-production}. 
Recognize that $u$ is monotonic along both time axes. Furthermore, the precise
path of $u$ minimizes the cost matrix at every step. Regions of increased 
cost in the cost matrix can be seen to repel the warp path. The 
distance $d$ between the demand and production curves here happens 
to be 0.756 GWe.

Dynamic time warping distance can therefore be used as an objective function 
to minimize for any demand and production curves. However, using full 
simulations to find $g(t, \Theta)$ remains expensive, even though DTW itself 
is computationally cheap. Therefore, a mechanism to reduce the overhead 
from production curve evaluation is needed.

\subsection{Gaussian Process Regression}
\label{gp}

Evaluating the production curve for a specific kind of facility using 
full fuel cycle simulations is relatively expensive, even in the 
computationally cheapest case of low-fidelity simulations. This is because a 
fuel cycle realization 
typically computes many features that, though coupled to the production 
curve, are not directly the production curve. For example, the mass balance of 
the fuel cycle physically bound the electricity production. However, the    
mass balances are not explicitly taken into account when trying to meet
a power demand curve.

Alternatively, surrogate models that predict the production curve directly
have many orders-of-magnitude fewer operations by virtue of not computing
implicit physical characteristics. This is not to say that the surrogate 
models are correct.  Rather, they are simply good enough to drive a demand
curve 
optimization. Surrogate models are used here inform a simulator about where
in the parameter space to look next. Truth about production curves should
still be derived from the fuel cycle simulator and not the surrogate model.
In the WORG algorithm, Gaussian processes are used to form the model. 

Gaussian processes are more fully covered elsewhere 
\cite{rasmussen2006gaussian}. Using Gaussian process for optimization has 
also been previously explored \cite{osborne2009gaussian}, though such studies 
tend not to 
investigate the integral problems posed by facility deployment. As with 
dynamic time warping, a minimal but sufficient introduction to GP regression
is presented 
for the purposes of the deployment optimization.
Consider the case of $Z$ simulations indexed by $z$ that each have a 
$\Theta_z$ deployment schedule and $g_z(t, \Theta_z)$ production curve.

A Gaussian process of these $Z$ simulations is set by its mean and 
covariance functions. The mean function is denoted as $\mu(t, \Theta)$ and 
is the expectation value $\E$ of 
the series of $G$ inputs:
\begin{equation}
\label{G}
G = \left\{g_1(t, \Theta_1), g_2(t, \Theta_2), \ldots, 
           g_Z(t, \Theta_Z)\right\}
\end{equation}
The covariance function is denoted $k(t, \Theta, t^\prime, \Theta^\prime)$ 
and is the expected value of the input to the mean. The mean and 
covariance can be expressed as
in Equations \ref{mean-func} \& \ref{covar-func} respectively.
\begin{equation}
\label{mean-func}
\mu(t, \Theta) = \E G
\end{equation}
\begin{equation}
\label{covar-func}
k(t, \Theta, t^\prime, \Theta^\prime) = 
    \E\left[(g_z(t, \Theta) - \mu(t, \Theta))
            (g_z(t^\prime, \Theta^\prime) - \mu(t^\prime, \Theta^\prime))
      \right]
\end{equation}
Note that in the above, the Gaussian process is itself $P+1$ dimensional, 
since the means and covariance are a function of both the deployment 
schedule ($P$) and time ($+1$).

The Gaussian process $\GP$ approximates the production curve 
given $Z$ simulations. Allow $*$ to indicate that the a quantity comes from 
the model as opposed to coming from the results of the simulator. A model 
production curve can then be written using either functional or operator
notation, as appropriate:
\begin{equation}
\label{gp-def-approx}
g_*(t, \Theta) \approx \GP\left(\mu(t, \Theta), 
                                 k(t, \Theta, t^\prime, \Theta^\prime)\right) 
                \equiv \GP G
\end{equation}
In machine learning terminology, $G$ serves as the training set for the 
GP model.

Now, when performing a regression on Gaussian processes, 
the nominal functional form for the covariance must be given. 
Such a functional form is also known as the the kernel function.
The kernel contains the \emph{hyperparameters} that are solved for to 
obtained a best-fit Gaussian process. The hyperparameters themselves are
defined based on the definition of the kernel function. Hyperparameter 
values are found via a regression of the maximal likelihood of 
the production curve. Any functional form could potentially serve as a kernel
function. However, a generally useful form is the is the exponential 
squared. This kernel can be seen in Equation \ref{exp2-kernel} with 
hyperparameters $\ell$ and $\sigma^2$ for a vector of parameters $r$:
\begin{equation}
\label{exp2-kernel}
k(r, r^\prime) = \sigma^2 \exp\left[-\frac{1}{2\ell}(r - r^\prime)^2 \right]
\end{equation}
However, other kernels such as the Mat\'ern $3/2$ kernel and Mat\'ern $5/2$
kernel \cite{paciorek2004nonstationary} were observed to be more robust for 
the WORG method. These can be seen in Equations \ref{matern-32} and 
\ref{matern-52} respectively.
\begin{equation}
\label{matern-32}
k(r, r^\prime) = \sigma^2 
                 \left(1 + \frac{\sqrt{3}}{\ell}|r - r^\prime|\right)
                 \exp\left(-\frac{\sqrt{3}}{\ell}|r - r^\prime|\right)
\end{equation}
\begin{equation}
\label{matern-52}
k(r, r^\prime) = \sigma^2 
                 \left(1 + \frac{\sqrt{5}}{\ell}|r - r^\prime|
                         + \frac{5}{3\ell^2}|r - r^\prime|^2\right)
                 \exp\left(-\frac{\sqrt{5}}{\ell}|r - r^\prime|\right)
\end{equation}

From here, say that $\K$ is a covariance matrix 
such that the element at the $r$-th row and $r^\prime$-th column is 
given by whichever kernel is chosen from 
Equations \ref{exp2-kernel}-\ref{matern-52}. Then the 
log likelihood $\log q$ of obtaining the training set production curves 
$G$ for a given time grid $\mathbf{t}$ and deployment schedule is as 
seen in Equation \ref{log-q}.
\begin{equation}
\label{log-q}
\log q(G|\mathbf{t}, \Theta) 
    = -\frac{1}{2}G^\top\left(\K + \tau^2\I\right)^{-1}G
      -\frac{1}{2}\log\left|\K + \tau^2\I\right|
      -\frac{ZTP}{2}\log 2\pi
\end{equation}
Here, $\tau$ is the uncertainty in the production curves coming from the 
simulations themselves. As most simulators do not report such uncertainties, 
$\tau$ may be set to floating point precision. $\I$ is the usual identity 
matrix. The hyperparameters $\ell$ and $\sigma^2$ are then adjusted via 
standard real-valued optimization methods such that Equation \ref{log-q} is 
as close to zero as possible. 
This regression of the Gaussian process itself yields the most likely 
model of the production curve knowing only a limited number of simulations.

However, the purpose of such a Gaussian process regression is to evaluate 
the production curve at points in time and for deployment schedules that 
have not been simulated. Take a time grid $\mathbf{t_*}$ and a hypothetical
deployment schedule $\Theta_*$. Now call the covariance vector between
the training set and the model evaluation    
$\mathbf{k}_* = \mathbf{k}(\mathbf{t_*}, \Theta_*)$. 
The production curve predicted by this Gaussian process is then given by
the following:
\begin{equation}
\label{metric-model}
\mathbf{g}_*(\mathbf{t}_*, \Theta_*) = 
    \mathbf{k}_*^\top \left(\K + \tau^2\I\right)^{-1}G
\end{equation}
Equations \ref{mean-func}-\ref{metric-model} are derived and discussed fully
in \cite{rasmussen2006gaussian}. 

Implementing the above Gaussian process mathematics for the specific
case of the WORG algorithm 
is not needed.  Free and open source Gaussian process modeling software 
libraries already exist and are applicable to the regression problem here.
Scikit-learn v0.17 \cite{scikit-learn} and George v0.2.1 \cite{hodlr} 
implement such a method and have a Python interface. George is specialized 
around Gaussian processes, and thus is preferred for WORG over scikit-learn, 
which is a general purpose machine learning library.

\begin{figure}[htb]
\centering
\includegraphics[width=0.9\textwidth]{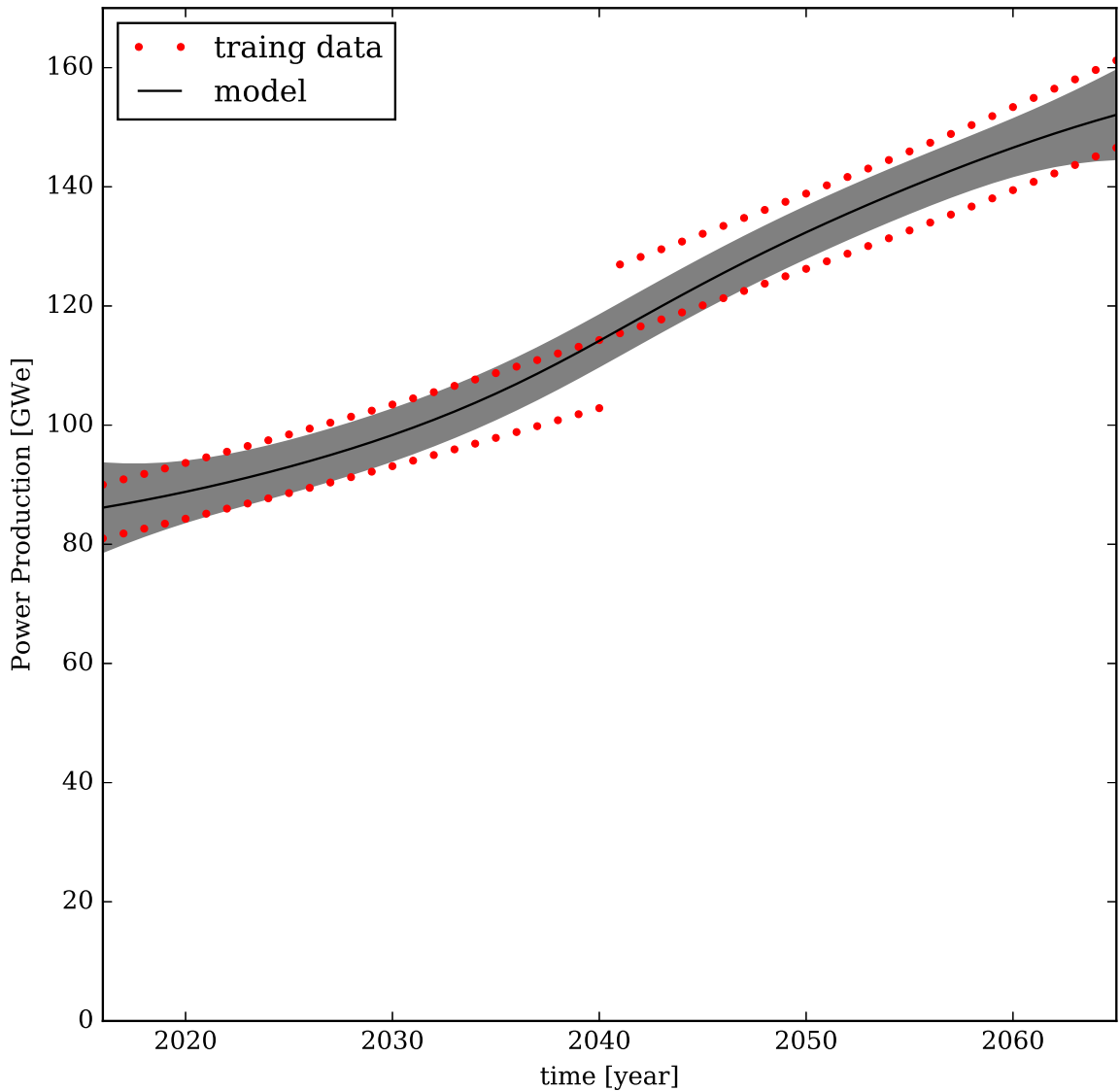}
\caption{The Gaussian process model of a 1\% growth curve along with the
an initial 10\% under production followed by a 10\% under production. 
The model is represented by the black line that runs between the red 
training points. Two standard deviations form the model are displayed as the
gray region.}
\label{gwe-model-}
\end{figure}

As an example, consider a Gaussian process between two power production 
curves similar to the example used in \S\ref{dtw}. The first is a nominal 1\% growth 
in GWe for 50 years starting at 
90 GWe in 2016. The second curve under-produces the first curve by 10\% 
for the first 25 years and over-produces by 10\% for the last 25 years.
Additionally, assume that there is a 10\% error on the training set data.
This will produce a model of the mean and covariance that splits the 
difference between these two curves. This example may be seen graphically
in Figure \ref{gwe-model-}.

The simple example above does not take advantage of an important 
feature of Gaussian processes. Namely, it is not limited to two production
curves in the training set.  As many as desirable may be used.  This will
allow the WORG algorithm to dynamically adjust the number of $Z$ simulations 
which are used to predict the next deployment schedule. WORG is thus capable
of 
effortlessly expanding $Z$ when new and useful simulations yield valuable production
curves.  However, it also enables $Z$ to contract to discard production
curves that would drive the deployment schedule away from an optimum.

Now that the Gaussian process regression and dynamic time warping tools have 
been added to the toolbox, the architecture of the WORG algorithm can 
itself be described.

\clearpage

\subsection{WORG Algorithm}
\label{algo}

The WORG algorithm has two fundamental phases on each iteration:
estimation and simulation.  These are preceded by an initialization 
before the optimization loop. Additionally, each iteration decides 
which information from the previous simulations is worth keeping for the
next estimation. Furthermore, the method of estimating deployment 
schedules may be altered each iteration.  Listing \ref{worg-pseudo}
shows the WORG algorithm as Python pseudo-code. A detailed walkthrough 
explanation of this code will now be presented.

Begin by initializing three empty sequences $\vec{\Theta}$, $G$, and $D$.
Each element of these series represents deployment schedule $\Theta$, 
a production curve $g(t, \Theta)$, and a dynamic time warping history 
between the demand and production curves $d(f, g)$.
Importantly, $\vec{\Theta}$, $G$, and $D$ only contain values for
the relevant optimization window $Z$. For example, root finding algorithms
such as Newton's method and the bisection method have a length-2 window
since they use the $(z-1)^\mathrm{th}$ point and the $z^\mathrm{th}$ point
to compute the $(z+1)^\mathrm{th}$ guess. Since a Gaussian process model is 
formed, any or all of the $s$ iterations may be used. However, 
restricting the optimization window to be either two or three depending on 
the 
circumstances balances the need to keep the points with the lowest $d$ 
values while pushing the model far from known regions with higher 
distances. Essentially, WORG tries to have $D$ contain one high-value $d$
and one or two low valued $d$ at all iterations. Such a tactic helps
form a meaningfully diverse GP model.

To this end, $\vec{\Theta}$, $G$, and $D$ are initialized with two 
bounding cases. The first is to set the deployment schedule equal to the
lower bound of the number of deployments $M$.  Recall that this is 
usually $\mathbf{0}$ everywhere, unless a minimum number of facilities 
must be deployed at a specific point in time. Running a simulation with 
$M$ will then yield a production curve $g(f, g)$ and the DTW distance to
this curve.  Note that just because the no facilities are deployed, the 
production curve need not be zero due to the initial conditions of the 
simulation. Existing initial facilities will continue to be productive. 

Similarly, another simulation may be executed for the maximum possible
deployment schedule $N$. This will also provide information on the 
production over time and the distance to the demand curve. $M$ and $N$
form the first two simulations, and therefore the loop 
variable $s$ is set to two.

\clearpage
\begin{lstlisting}[
    caption={WORG Algorithm in Python Pseudo-code},
    label=worg-pseudo,mathescape]
Thetas, G, D = [], [], []  # initialize history

# run lower bound simulation
g, d = run_sim(M, f)
Thetas.append(M)
G.append(g)
D.append(d)

# run upper bound simulation    
g, d = run_sim(N, f)
Thetas.append(N)
G.append(g)
D.append(d)

s = 2
while MAX_D < D[-1] and s < S:
    # set estimation method
    method = initial_method
    if method == 'all' and (s%4 < 2):
        method = 'stochastic'

    # estimate deployment schedule and run simulation
    Theta = estimate(Thetas, G, D, f, method)
    g, d = run_sim(Theta, f)
    Thetas.append(Theta)
    G.append(g)
    D.append(d)

    # take only the most important and most recent schedules
    idx = argsort(D)[:2]
    if D[-1] == max(D):
        idx.append(-1)
    Thetas = [Thetas[i] for i in idx]
    G = [G[i] for i in idx]
    D = [D[i] for i in idx]
    s = (s + 1)
\end{lstlisting}
\clearpage

The optimization loop may now be entered.  This loop has two conditions.
The first is that the next iteration occurs only if the last distance
is greater than a threshold value $\mathrm{MAX\_D}$. The second is that 
the loop variable $s$ must be less than the maximum number
of iteration $S$.

The first step in each iteration is to choose the estimation method. The
three mechanisms will be discussed in detail in \S\ref{selecting}. For
the purposes of the optimization loop, they may be represented by the 
\stochastic, \innerprod, and \allflag flags. The stochastic method 
chooses many random deployment schedules to test. Alternatively, an inner
product search of the space defined by $M$ and $N$ may be performed. 
Lastly, the \allflag flag performs both of the previous estimates and takes
the one with lowest computed distance.  However, \allflag can sometimes
declare the inner product search the winner for all $s$.  This can itself 
be problematic since this estimation method has the tendency to form 
deterministic loops when close to an optimum. This behavior is not unlike
similar loops formed with floating point approximations to Newton's method.
To prevent this when using \allflag, WORG forces the stochastic method
for two consecutive iterations out of every four.  
 
A best-guess estimate for a deployment schedule $\Theta$ may finally be
made.  This takes the previous deployment schedules $\vec{\Theta}$ and
production curves $G$ and forms a Gaussian process model. Potential 
values for $\Theta$ are explored according to the selected estimation method.  The 
$\Theta$ that produces the minimum dynamic time warping distance between
the demand curve and the model $d(f, g_*)$ is then returned.

The $\Theta$ estimate is then supplied to the the simulator itself and 
a simulation is executed.  The details of this procedure are, of course,
simulator specific.  However, the simulation combined with any post-processing 
needed should return an aggregate production curve $g_s(t, \Theta)$.  
This is then compared to demand curve via $d(f, g_s)$. After the simulation, 
$\Theta$, $g_s(t, \Theta)$, and $d(f, g_s)$ are appended to the 
$\vec{\Theta}$, $G$, and $D$ sequences.  Note that the 
production curve and DTW distance from the simulator are appended, 
not the production curve and distance from the model estimate.

Concluding the optimization loop, $\vec{\Theta}$, $G$, $D$, and $s$ are
updated.  This begins by finding and keeping the two elements with the 
lowest distances between the demand and production curves.  However, 
if the most recent simulation yielded the largest distance, this is also
kept for the next iteration. Keeping the largest distance serves to deter
exploration in this direction on the next iteration.  Thus a sequence of 
two or three indices is chosen. These indices are applied to redefine
$\vec{\Theta}$, $G$, and $D$. Lastly, $s$ is incremented by one and the
next iteration begins.

The WORG algorithm presented here shows the overall structure of the 
optimization.  However, equally important and not covered in this section 
is how the
estimation phase chooses $\Theta$.  The methods that WORG may use are 
presented in the following section and completes the methodology.

\section{Selecting Deployment Schedule Estimates}
\label{selecting}

There are three methods for choosing a new deployment schedule $\Theta$ 
to attempt to run in a simulator. The first is stochastic with weighted
probabilities for the $\theta_p$.  The second does a deterministic sweep 
iteratively over all options, minimizing the dynamic time warping distance
at each point in time for each deployment parameter.  The last combines
these two and choose the one with the minimum distance to the demand curve.

All of these rely on a Gaussian process model of the production
curve. This is because constructing and evaluating GP model $g_*$ is 
significantly faster than performing even a low-fidelity simulation. 
As a demonstrative example, say each evaluation of $d(f, g_*)$ takes a tenth
of a second (which is excessively long) and $d(f, g_s)$ for a low fidelity
simulation takes ten seconds (which is reasonable), the model evaluation 
is still one hundred times faster.  Furthermore, the cost of constructing 
the GP model can is amortized over the number of guesses that are made.

However, the choice of which $\theta_p$ to pick is extremely important
as they drive the optimization. In a vanilla stochastic algorithm, 
each $\theta_p$ would be selected as a univariate integer on the 
range $[M_p, N_p]$.  However, this ignores the distance information $D$ 
that is known about the training set which is used to create the Gaussian 
process.
More intelligent guesses for $\theta_p$ focus the model 
evaluations to more promising regions of the option space.  This in turn 
helps reduce the overall number of expensive simulations needed to find 
a `good enough' deployment schedule.

The three WORG $\Theta$ selection methods are described in order in the 
following subsections.

\subsection{Stochastic Estimation}
\label{stochastic}

The stochastic method works by randomly choosing $\Gamma$ deployment 
schedules and evaluating $g_*(t, \Theta_\gamma)$ for each guess $\gamma$.
The $\Theta_\gamma$ which has the minimum distance $d_\gamma$ is 
taken as the best-guess deployment schedule.  The number of guesses may 
be as large or as small as desired.  However, a reasonable number to pick 
spans the option space. This is simply is the $L_1$ norm of the difference 
of the inclusive bounds. Namely, set $\Gamma$ as in Equation 
\ref{Gamma-default} for a minimum number of for stochastic guesses.
\begin{equation}
\label{Gamma-default}
\Gamma = \sum_p^P (N_p - M_p + 1)
\end{equation}
Each $\theta_p$ has $N_p - M_p + 1$ options. Thus a reasonable
choice for $\Gamma$ is the sum of the number of independent options.

Still, each option for $\theta_p$ should not be equally likely. 
For example, if the demand curve is relatively low, the number of deployed 
facilities is unlikely to be relatively high. For this reason, the choice 
of $\theta_p$ should be weighted.  Furthermore, note that each $\theta_p$
is potentially weighted differently as they are all independent parameters.
Denote $n \in [M_p, N_p]$ such that the n-th weight for the p-th parameter 
is called $w_{n,p}$. 

To choose weights, first observe that the distances $D$ can be said to be
inversely proportional to how likely each deployment schedule in 
$\vec{\Theta}$ should be. A one-dimensional Gaussian process can thus be
constructed to model inverse distances given the values of the deployment 
parameter for each schedule, namely $\vec{\theta_p}$.  Call this model 
$d_*^{-1}$ as seen in Equation \ref{d-inv-model}.
\begin{equation}
\label{d-inv-model}
d_*^{-1}(\theta_p) = \GP\left(\mu(\vec{\theta_p}), 
                              k(\vec{\theta_p}, \vec{\theta_p}^\prime)\right)
                   \equiv \GP\left[D^{-1}\right]
\end{equation}
The construction, regression of hyperparameters, and evaluation of this 
model follows analogously to the production curve modeling presented in 
\S\ref{gp}.

The weights for $\theta_p$ are then the normalized evaluation of the 
inverse distance model for all $m$ and $n$ defined on the p-th range.
Symbolically, 
\begin{equation}
\label{d-inv-w}
w_{n,p} = \frac{d_*^{-1}(n)}{\sum_{m=M_p}^{N_p} d_*^{-1}(m)}
\end{equation}
Equation \ref{d-inv-w} works very well as long as a valid model can 
be established.  However, this is sometimes not the case when the $\theta_p$
are degenerate, the distances are too close together, the distances are 
too close to zero, or other stability issues arise.

In cases where a valid model may not be formed for $d_*^{-1}(\theta_p)$, 
a Poisson distribution may be used instead.  Take the mean of the Poisson
distribution $\lambda$ to be the value of $\theta_p$ where the distance
is minimized.
\begin{equation}
\label{lambda}
\lambda_p = \theta_p | \mathrm{argmin}(D) 
\end{equation}
Hence, the Poisson probability distribution for the $n$-th weight of the
$p$-th deployment parameter is, 
\begin{equation}
\label{poisson}
\mathrm{Poisson(n)} = \frac{(\lambda_p)^n}{n!} e^{-\lambda_p}
\end{equation}
Now, because $n$ is bounded, it is important to renormalize Equation 
\ref{poisson} when constructing stochastic weights.
\begin{equation}
\label{poisson-w}
\begin{split}
w_{n,p} & = \frac{\frac{(\lambda_p)^n}{n!} e^{-\lambda_p}}
                 {\sum_{m=M_p}^{N_p} \frac{(\lambda_p)^m}{m!} e^{-\lambda_p}}\\
        & = \frac{(\lambda_p)^n}
                 {n!\sum_{m=M_p}^{N_p} \frac{(\lambda_p)^m}{m!}}\\
\end{split}
\end{equation}
Poisson-based weights could be used exclusively, foregoing the inverse 
distance Gaussian process models completely. However, a Poisson-only 
method takes into account less information about the demand-to-production
curve distances. It was therefore observed to converge more slowly 
on an optimum than using Poisson weights as a backup.  Since the total 
number of simulations is aiming to be minimized for \emph{in situ} use, 
the WORG method uses Poisson weights as a fallback only.

After weights are computed for all $P$ deployment parameters, a set of 
$\Gamma$ deployment schedules may be stochastically chosen. The Gaussian
process for each $g_*(\mathbf{t}, \Theta_\gamma)$ is then evaluated and the
dynamic time warping distance to the demand curve is computed. The 
deployment schedule with the minimum distance is then selected and returned.

\subsection{Inner Product Estimation}
\label{inner-prod}

As an alternative to the stochastic method demonstrated in \S\ref{stochastic}, 
a best-guess for $\Theta$ can also be built up iteratively over all times.
The method here uses the same production curve Gaussian process $g_*$ to 
predict production levels and measure the distance to the demand curve.
However, this method minimizes the distance at time
$t$ and then uses this to inform the minimization of $t+1$. Starting at $t=1$
and moving through the whole time grid to $t=T$, a complete deployment 
schedule is generated.

The following description is for the simplified case when $P==T$. However, 
this method is easily extended to the case where $P > T$, such as for 
multiple reactor types.  When $P > T$, group $\theta_p$ that occur on the 
same time step together and take the outer product of their options prior to
stepping through time.

For this method, define the time sub-grid $\mathbf{t}_p$ as the sequence of
all times less than or equal to the time when parameter $p$ occurs, $t(p)$.
\begin{equation}
\label{t-p}
\mathbf{t}_p = \left\{t | t \le t(p)\right\}
\end{equation}
Now define the deployment schedule $\Theta^t$ up to time $t$ through the
following recursion relations:
\begin{equation}
\label{Theta-t}
\begin{split}
\theta_1 & = n \, | \, \mathrm{min}\left[d(f, g_*(1, n))\right]
                       \forall n\in[M_1, N_1] \\
\Theta^1 & = \left\{\theta_1\right\}\\
\theta_p & = n \, | \, \mathrm{min}\left[d(f, g_*(\mathbf{t}_p, 
                                                  \Theta^{t-1}, n))\right]
                       \forall n\in[M_p, N_p] \\
\Theta^t & = \left\{\Theta_1^{t-1}, \ldots, \Theta_{p-1}^{t-1}, \theta_p\right\}
\end{split}
\end{equation}
Equation \ref{Theta-t} has the effect of choosing the the number of 
facilities to deploy at each time step that minimizes the distance function.
The current time step uses the previous deployment schedule and only searches
the option space of the its own deployment parameter $\theta_p$. 
Once $\Theta^T$ is 
reached, it is selected as the deployment schedule $\Theta$. The inner 
product method here requires the same number of model evaluations of $g_*$ as 
were selected for the default value of $\Gamma$ in Equation 
\ref{Gamma-default} for stochastic estimation.

\subsection{All of the Above Estimation}
\label{all}

This method is simply to run both the stochastic method and the inner product
method and determine which has the lower $d(f, g_*)$ for the deployment 
schedules they produce.  This method contains both the advantages and 
disadvantages of its constituents.  Additionally, it has the disadvantage 
of being more computationally expensive than the other methods individually.

The advantage from the stochastic method is that the entire space is 
potentially searched. There are no forbidden regions.
This is important since there may be other
optima far away from the current $\vec{\Theta}$ that produce lower distances.
Searching globally prevents the stochastic method from becoming stuck locally.
However, the stochastic method may take many iterations to 
make minor improvements on a $\Theta$ which is already close to a best-guess.
It is, after all, searching globally for something better.

On the other hand, the inner product method is designed to search around 
the part of the Gaussian process model which already produces good results. It
is meant to make minor adjustments as it goes.  Unfortunately, this means
the inner product method can more easily get stuck in a cycle  where it 
produces the same series of deployment schedules over and over again. 
It has no mechanism on its own to break out of such cycles.

With the all-of-the-above option, the job of balancing the relative merits
of the stochastic and inner product methods is left to the optimization 
loop itself.  This can be seen in \S\ref{algo}.  If the \allflag flag 
is set as the estimation method, it is only executed as the \allflag flag 
two of every four 
iterations.  Other strategies for determining how and when each of the 
three methods are used could be designed. However, any more complex strategy 
should be able to show that it meaningfully reduces the number of 
optimization loop iterations required.

At this point, the entire WORG method has been described. A
demonstration of how it performs for a representative fuel cycle is 
presented in the next section. 

\section{Results \& Performance}
\label{results}

To demonstrate the three variant WORG methods, an unconstrained 
once-through fuel cycle is modeled with the Cyclus simulator 
\cite{DBLP:journals/corr/HuffGCFMOSSW15}. In such a scenario, uranium
mining, enrichment, fuel fabrication, and storage all have effectively 
infinite capacities. The only meaningful constraints on the system are
how many light-water reactors (LWR) are built.

The base simulation begins with 100 reactors in 2016 that each produce
1 GWe, have an 18 month batch length with a one month reload time.
The initial fleet of LWRs retires evenly over the 40 years from 2016 to 
2056. All new reactors have 60 year life times.  The simulation itself 
follows 20 years from 2016 to 2035. This is on the higher end of 
\emph{in situ} time horizons expected, which presumably 
will be in the neighborhood of 1, 5, 10, or 20 years.

The study here compares how WORG performs for 0\% (steady state), 1\%, 
and 2\% growth curves from an initial 90 GWe target. These are examined
using the three estimation methods variants described in the previous section.
Calling $\rho$ the growth rate as a 
fraction, the demand curve is thus,
\begin{equation}
\label{f-rate}
f(t) = 90 (1 + \rho)^t
\end{equation}
Moreover, the upper bound for the number of deployable facilities at 
each time is set to be the ceiling of ten times the total growth. 
That is, assuming ten facilities at most could be deployed in the first
year, increase the upper bound along with the growth rate.  This yields
the following expression for $N$.
\begin{equation}
\label{n-rate}
N(t) = \left\lceil 10 (1 + \rho)^t\right\rceil
\end{equation}
The lower bound for the number of deployed reactors is taken to be the 
zero vector, $M = \mathbf{0}$.  A maximum of twenty simulations are allowed, 
or $S = 20$.
This is because an \emph{in situ} method cannot afford many optimization 
iterations. The random seed for all optimizations was 424242.

Note that because of the integral nature of facility deployment, 
exactly matching a continuous demand curve is not possible in general. 
Slight over- and under-prediction are expected for most points in time. 
Furthermore, 
it is unlikely that the initial facilities will match the demand curve 
themselves. If the initial facilities do no meet the demand on their own, 
then the optimized deployment schedule is capable making up the difference.
However if the initial facilities produce more than the demand curve, 
the the optimizer is only capable of deploying zero facilities at these
times. The WORG method does not help make radical adjustments to 
accommodate 
problems with the initial conditions, such as when 50 GWe are demanded 
but 100 GWe are already being produced.

\begin{figure}[htb]
\centering
\includegraphics[width=0.9\textwidth]{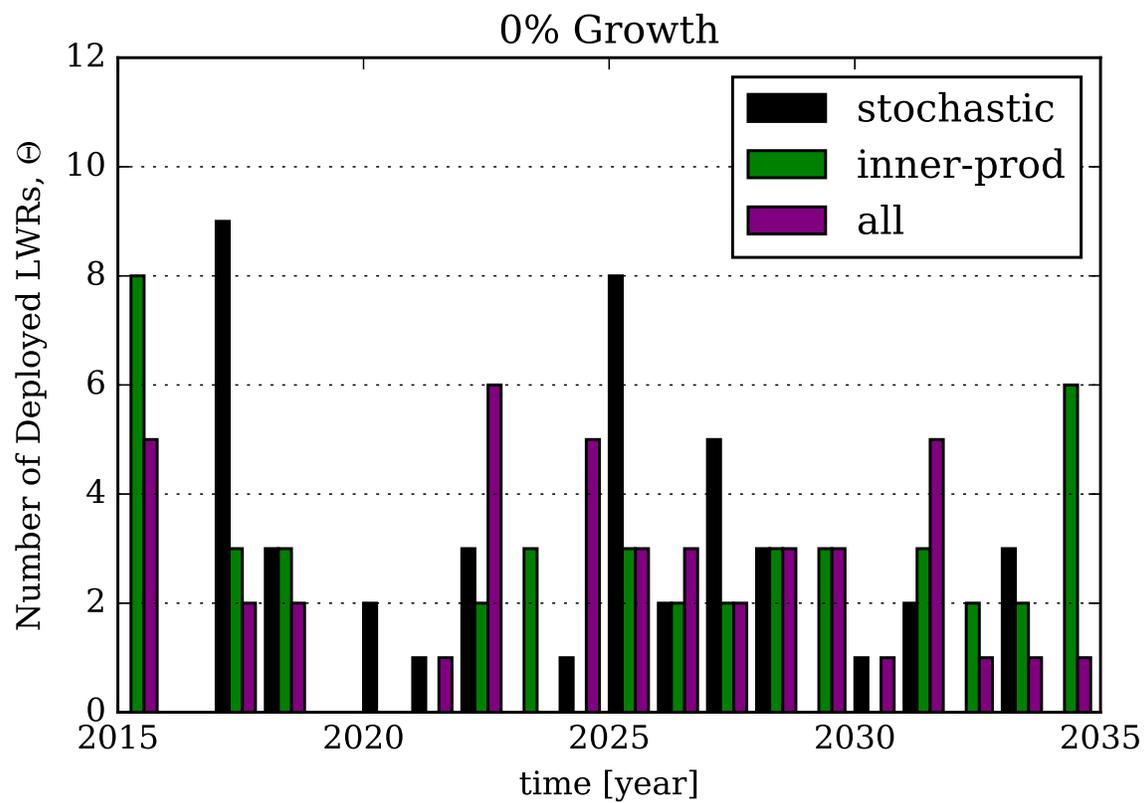}
\caption{Optimized deployment schedule $\Theta$ for a 0\% growth (steady 
state) demand curve. The number of deployed facilities shown are 
for the \stochastic (black), \innerprod (green), and \allflag (purple)
estimation methods.}
\label{deploy-0}
\end{figure}

\begin{figure}[htb]
\centering
\includegraphics[width=0.9\textwidth]{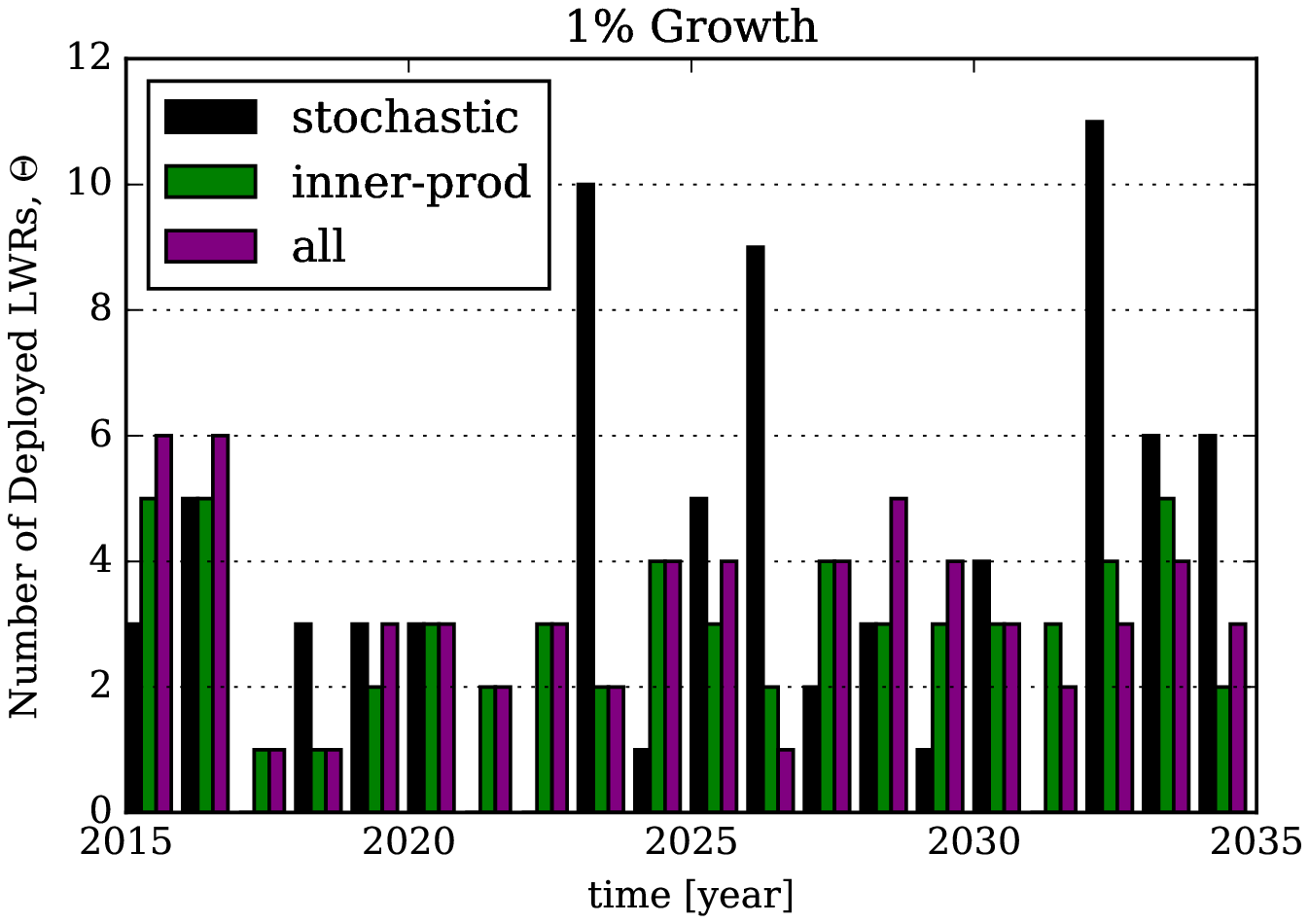}
\caption{Optimized deployment schedule $\Theta$ for a 1\% growth 
demand curve. The number of deployed facilities shown are 
for the \stochastic (black), \innerprod (green), and \allflag (purple)
estimation methods.}
\label{deploy-1}
\end{figure}

\begin{figure}[htb]
\centering
\includegraphics[width=0.9\textwidth]{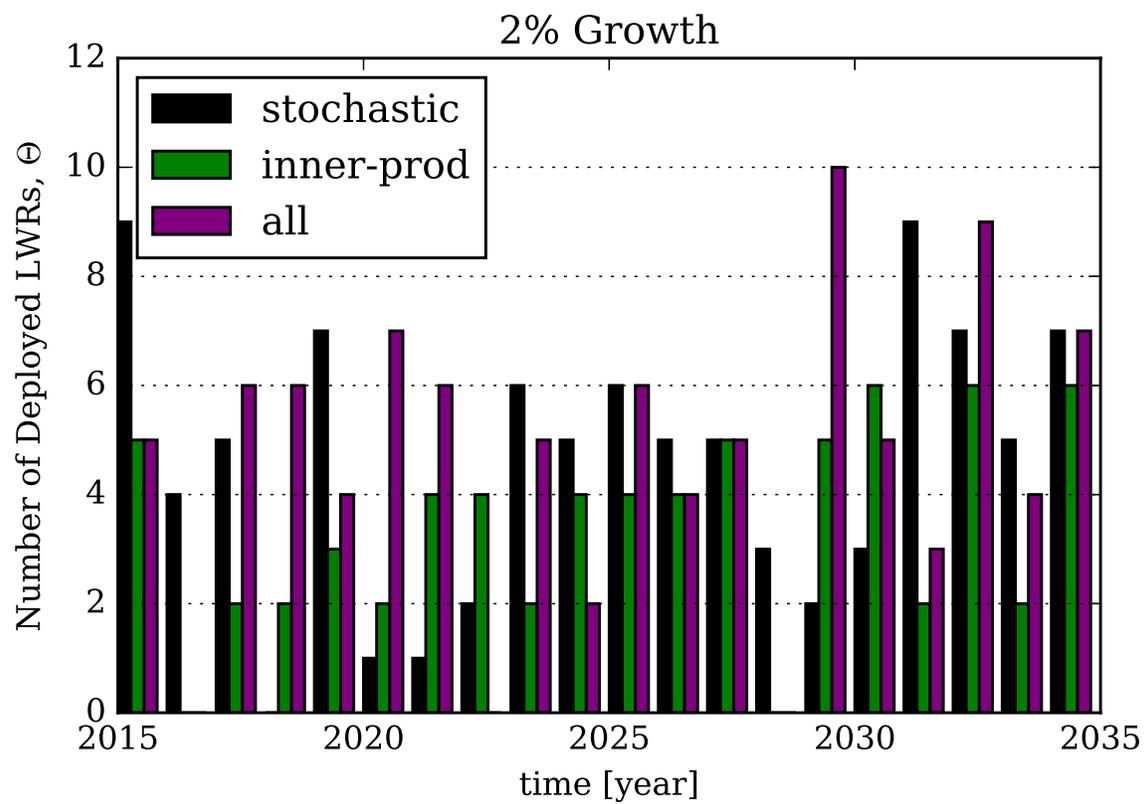}
\caption{Optimized deployment schedule $\Theta$ for a 2\% growth 
demand curve. The number of deployed facilities shown are 
for the \stochastic (black), \innerprod (green), and \allflag (purple)
estimation methods.}
\label{deploy-2}
\end{figure}

First examine Figures \ref{deploy-0} - \ref{deploy-2} which
show the optimized deployment schedule for the three estimation methods.
The figures represent the 0\%, 1\% and 2\% growth curves respectively.
Figures \ref{deploy-0} \& \ref{deploy-1} show that the \stochastic method
has the highest number of facilities deployed at any single time. In the 
2\% case, the \allflag method predicts the largest number of facilities 
deployed.

\clearpage

\begin{figure}[htb]
\centering
\includegraphics[width=0.9\textwidth]{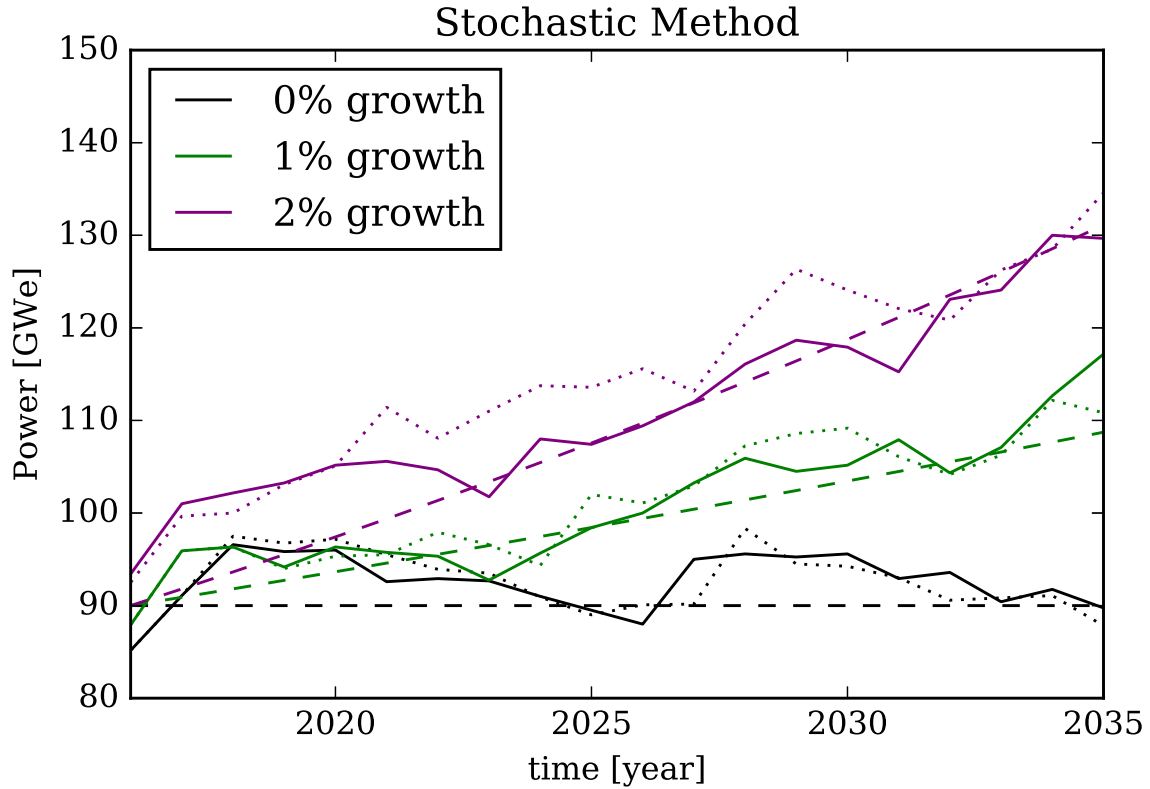}
\caption{Power production to demand comparison for 20 year deployment 
schedule optimization using only the \stochastic estimation method.
0\%, 1\%, and 2\% growth rates starting at 90 GWe are shown. Solid lines 
represent the best guess deployment schedule.  Dotted lines are represent 
the second best guess deployment schedule. Dashed lines represent the 
demand curve that is targeted.
}
\label{demand-product-stochastic}
\end{figure}

\begin{figure}[htb]
\centering
\includegraphics[width=0.9\textwidth]{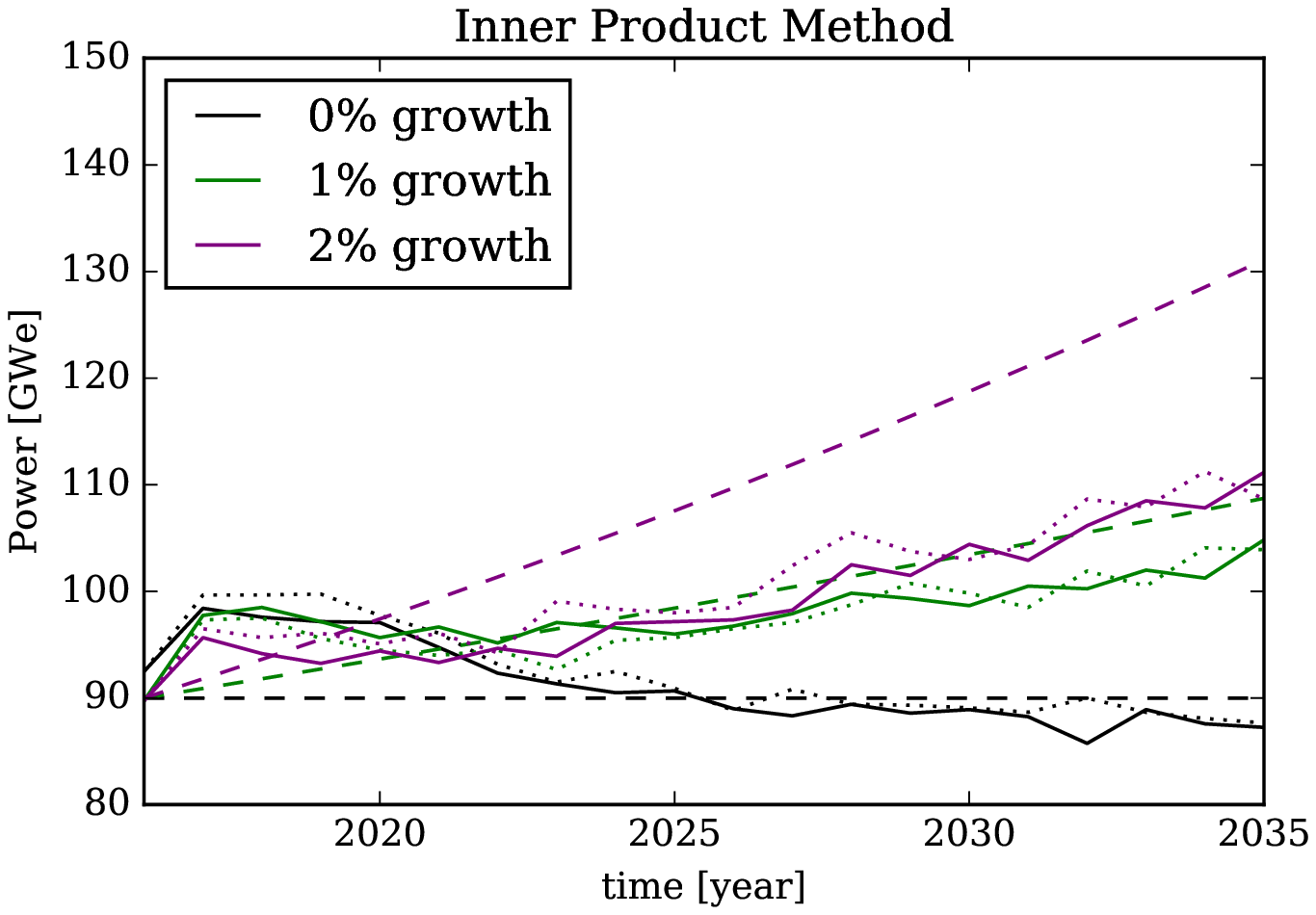}
\caption{Power production to demand comparison for 20 year deployment 
schedule optimization using only the \innerprod estimation method.
0\%, 1\%, and 2\% growth rates starting at 90 GWe are shown. Solid lines 
represent the best guess deployment schedule.  Dotted lines are represent 
the second best guess deployment schedule. Dashed lines represent the 
demand curve that is targeted.
}
\label{demand-product-inner-product}
\end{figure}

\begin{figure}[htb]
\centering
\includegraphics[width=0.9\textwidth]{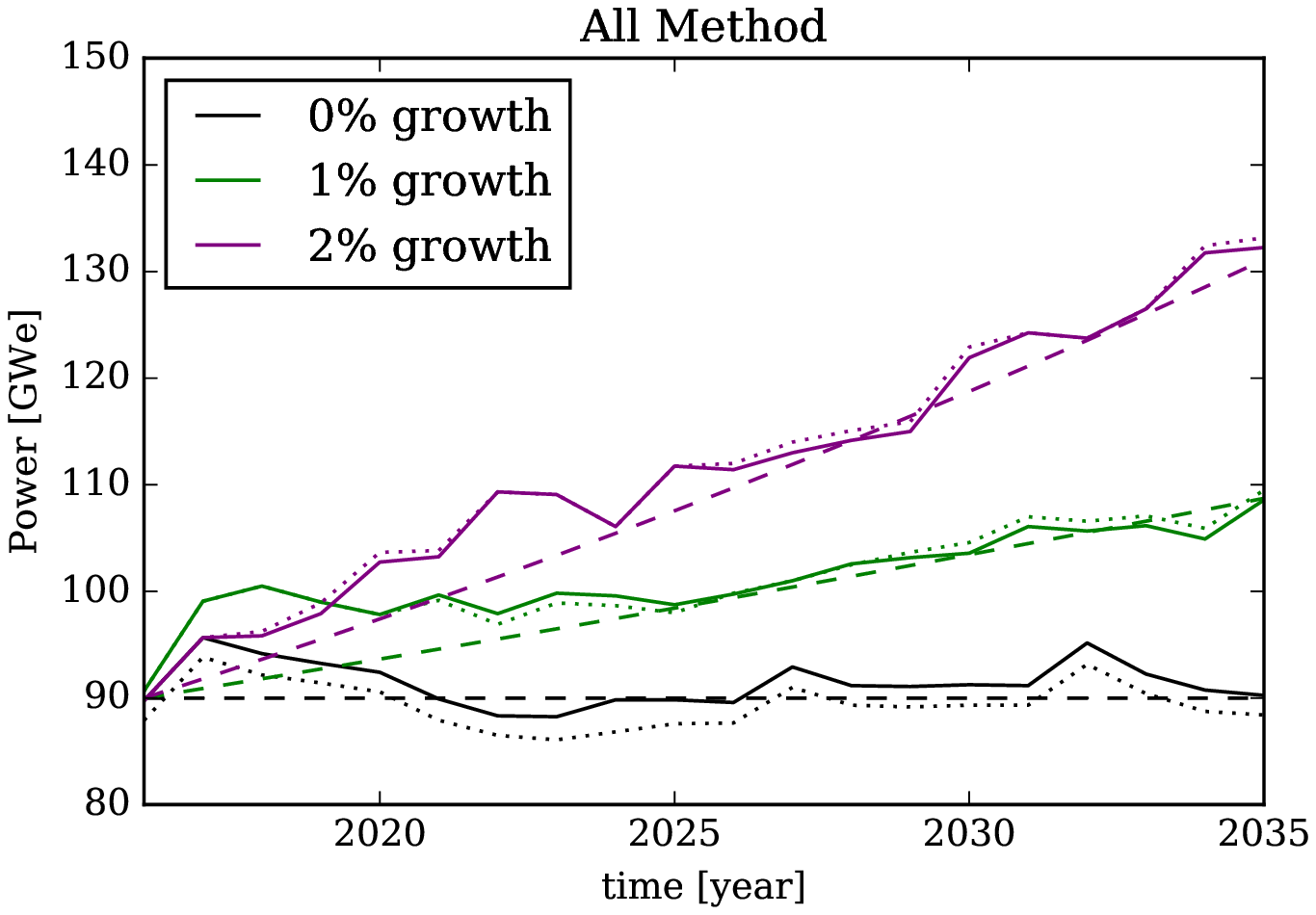}
\caption{Power production to demand comparison for 20 year deployment 
schedule optimization using the \allflag estimation method that selects the 
best of both \stochastic and \innerprod estimations.
0\%, 1\%, and 2\% growth rates starting at 90 GWe are shown. Solid lines 
represent the best guess deployment schedule.  Dotted lines are represent 
the second best guess deployment schedule. Dashed lines represent the 
demand curve that is targeted.}
\label{demand-product-all}
\end{figure}

More important than the deployment schedules themselves, however, are the
production curves that they elicit.
Figures \ref{demand-product-stochastic} - 
\ref{demand-product-all} display the
power production for the best-guess deployment schedule $G_1$ (solid lines),
production for the second best schedule $G_2$ (dotted lines), 
and demand curves (dashed lines) for 0\%, 1\%, and 2\% growth. 
The figures show each estimation mechanism separately.

As seen in Figure \ref{demand-product-stochastic}, the \stochastic only 
estimations follow the trend line of the growth curve.  However, only for
a few regions such as for 2\% growth between 2025 - 2030, does the 
production very closely match the the demand.  The second best guess for the 
deployment schedule shows a relatively large degree of eccentricity. This 
indicates that $G_2$ at the end of the optimization is there to show the 
Gaussian process model regions in the $\Theta$ option space that do not work. 

Alternatively the inner product search estimation can be seen in Figure 
\ref{demand-product-inner-product}. This method does a reasonable job of
predicting a steady state scenario at later times. However, the 1\% and 2\%
curves are largely under-predicted.  For the 2\% case, this is so severe
as to be considered wholly wrong. This situation arises because the 
\innerprod method can enter into deterministic traps where the same cycle
of deployment schedules is predicted forever.  If such a trap is fallen
into when the 
production curve is far from an optimum, the inner product method does not
yield a sufficiently close guess of the deployment schedule.

However, combining the stochastic and inner product estimation methods limits
the weakness of each method. Figure \ref{demand-product-all} shows the 
production and demand information for the \allflag estimation flag. 
By inspection, 
this method produces production curves that match the demand much more 
closely. This is especially true for later times. Over-prediction 
discrepancies for early times come from the fact the initially deployed 
facilities (100 LWRs with 18 month cycles) does not precisely match an 
initial 90 GWe target. The problem was specified this way in order to show 
that reasonable deployments are still selected even in slightly 
unreasonable situations.

Still, the main purpose of the WORG algorithm is to converge as quickly as 
possible to a reasonable best-guess $\Theta$. The limiting factor is the 
number of 
predictive simulations which must be run. WORG will always execute 
at least three full simulations: the lower bound, the upper bound, and one
iteration of the optimization loop. A reasonable limit on the total number 
of simulations $S$ is 20.  For an \emph{in situ} calculator, it is unlikely
to want to compute more than this number of sub-simulations to predict 
the deployment schedule for the next 1, 5, 10, or 20 years. However, it 
may be possible to limit $S$ to significantly below 20 as well.

\begin{figure}[htb]
\centering
\includegraphics[width=0.9\textwidth]{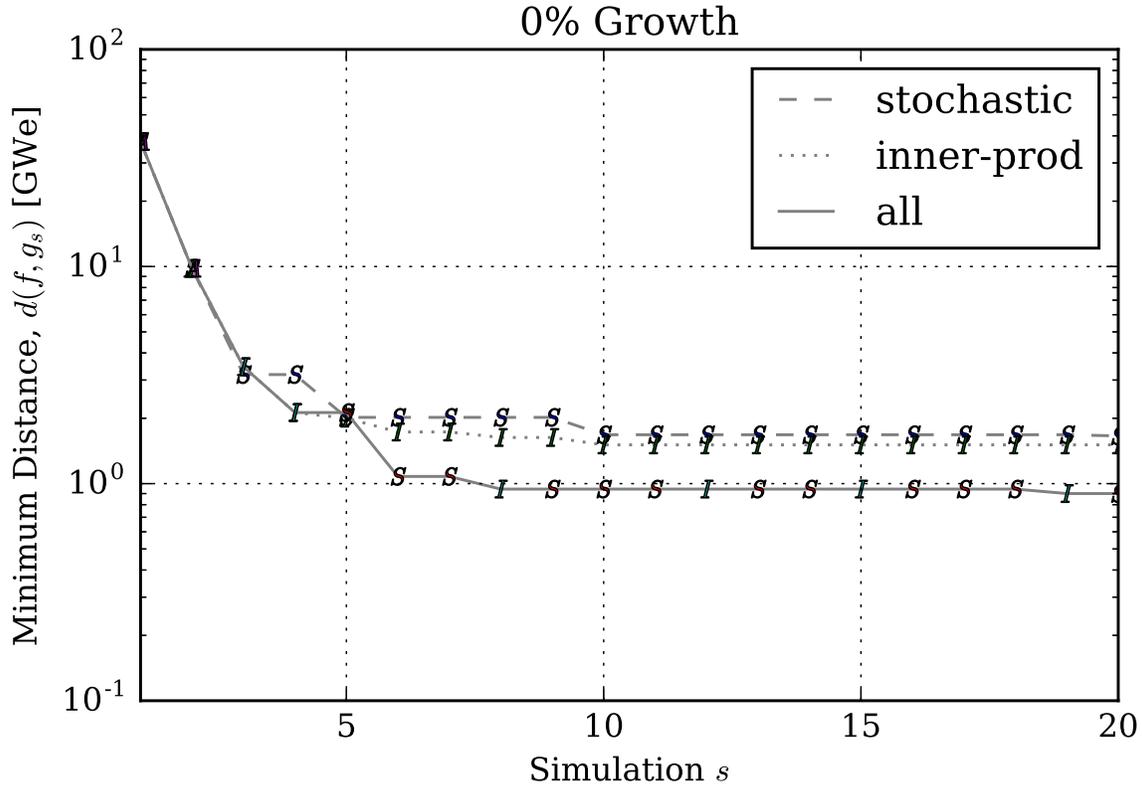}
\caption{Convergence of 0\% growth rate solution for the distance between
the demand and production curves $d(f, g)$ as a function of the number of 
simulations. The three estimation methods are shown. Additionally, the
$S$, $I$, and $A$ marker represent whether the \stochastic, \innerprod, 
or \allflag method was selected as the best fit. For $2 < s$, the \allflag
will select either $S$ or $I$.
}
\label{converge-0per}
\end{figure}

\begin{figure}[htb]
\centering
\includegraphics[width=0.9\textwidth]{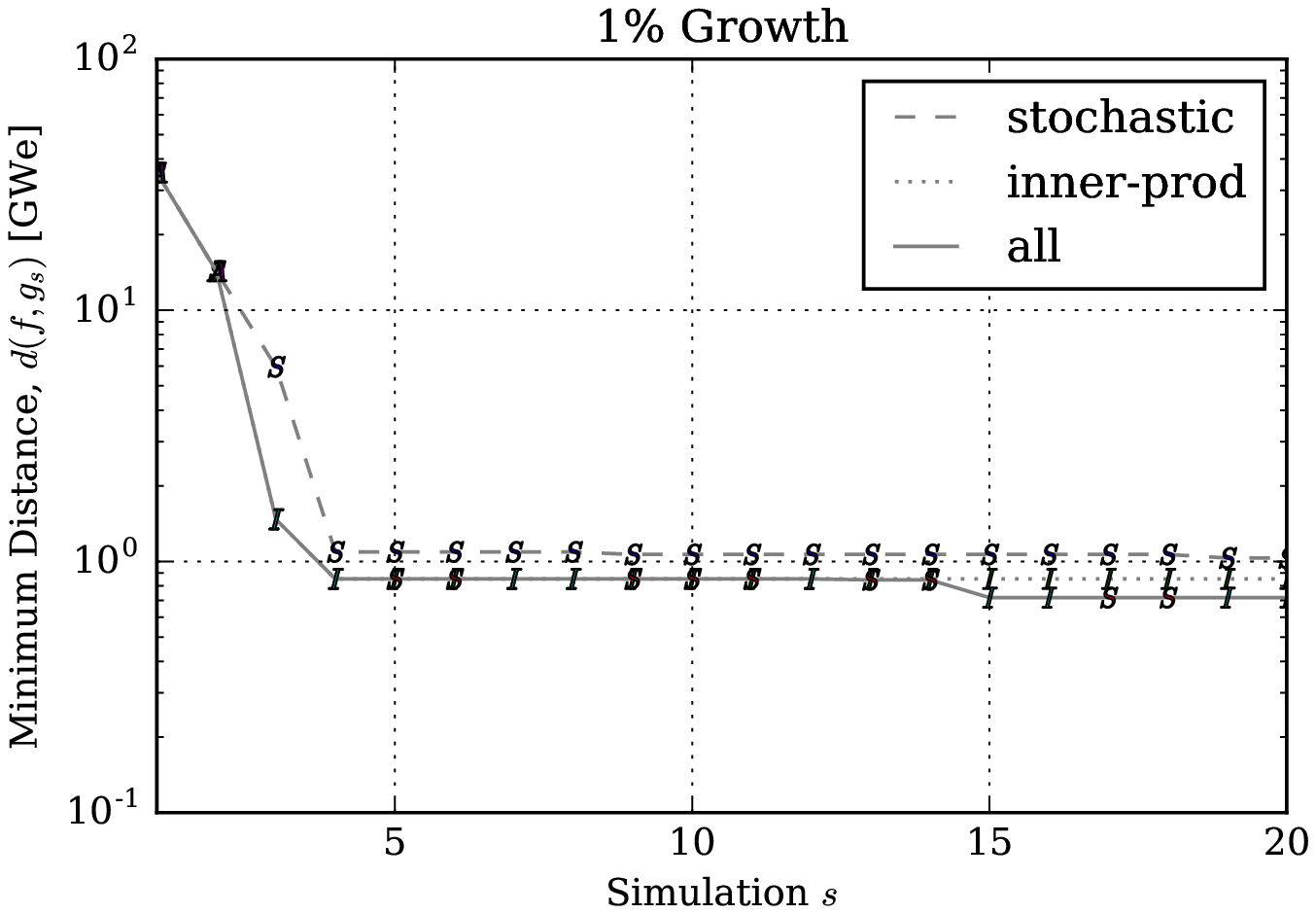}
\caption{Convergence of 1\% growth rate solution for the distance between
the demand and production curves $d(f, g)$ as a function of the number of 
simulations. The three estimation methods are shown. Additionally, the
$S$, $I$, and $A$ marker represent whether the \stochastic, \innerprod, 
or \allflag method was selected as the best fit. For $2 < s$, the \allflag
will select either $S$ or $I$.
}
\label{converge-1per}
\end{figure}

\begin{figure}[htb]
\centering
\includegraphics[width=0.9\textwidth]{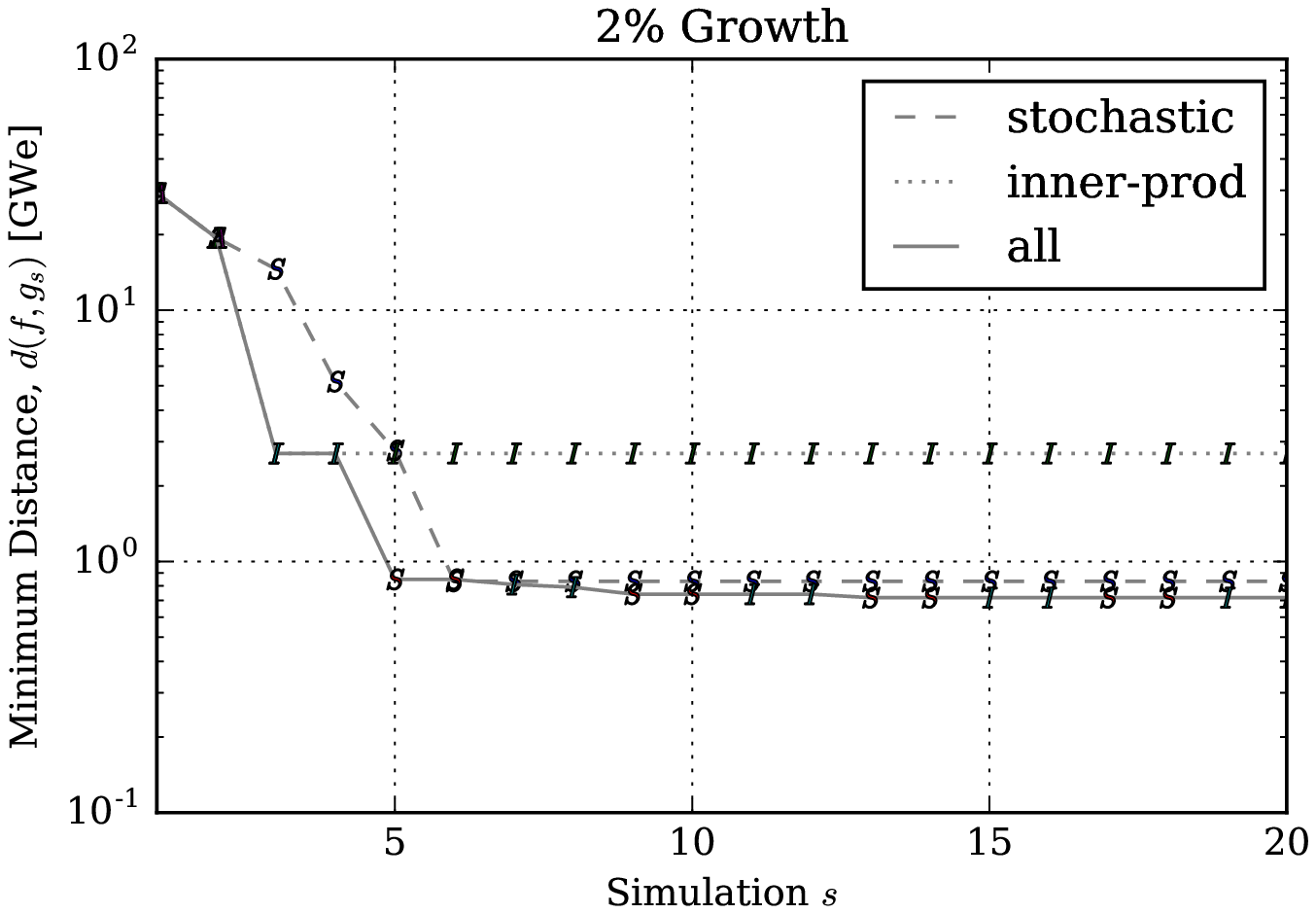}
\caption{Convergence of 2\% growth rate solution for the distance between
the demand and production curves $d(f, g)$ as a function of the number of 
simulations. The three estimation methods are shown. Additionally, the
$S$, $I$, and $A$ marker represent whether the \stochastic, \innerprod, 
or \allflag method was selected as the best fit. For $2 < s$, the \allflag
will select either $S$ or $I$.
}
\label{converge-2per}
\end{figure}

Figures \ref{converge-0per} - \ref{converge-2per} demonstrate a convergence 
study for the three 
different growth rates as a function of the number of simulations $s$.
These three figures demonstrate important properties of the WORG method.
The first is that the \allflag estimation method here consistently has
the lowest minimum dynamic time warping distance, and thus the best guess 
for the deployment schedule $\Theta$.  Furthermore, the \allflag curves
in these figures show that the best estimation method consistently switches
between the inner product search and stochastic search. These imply 
that together the \stochastic and \innerprod methods converge more quickly 
that the sum of their parts.  This is particularly visible in the 0\% growth
rate case seen in Figure \ref{converge-0per}.

Additionally, These convergence plots show that the majority of the 
$\Theta$ selection gains are made by $s=10$.  Simulations past ten may 
show differential improvement for the \allflag method.  However, they 
do not tend to generate any measurable improvement for \stochastic or
\innerprod selection mechanisms. Furthermore, most of the gains for the 
\allflag method are realized by simulation five or six. Thus a reduction
by a factor of two to four in the number of simulations is available.
This equates directly to a like reduction in computational cost.

Moreover,
Figures \ref{converge-0per} - \ref{converge-2per} also show the tendency of
the \innerprod method to become deterministically stuck when used on its own.
In all three cases the \innerprod method resolves to a constant by 
simulation five or ten.  In the 1\% case, this solution happens to be quite close
to the solution predicted by the \allflag method.  In the 0\% case, this
constant is meaningfully distinct, but is still more akin to the \stochastic 
solution than the \allflag solution.  In the 2\% case, the converged \innerprod
result has little to do with the \allflag prediction.  Thus it is not
recommended to use \innerprod on its own. While it may succeed, it is too 
risky because it may cease improving at the wrong place.

The \stochastic estimation method is also prone to similar cyclic behavior.
However, in the \innerprod method, 
is the second best guess $G_2$ is also in a constant prediction cycle along 
with $G_1$.  With the stochastic search, the second best 
guess has more freedom to roam the option space, creating potentially 
different Gaussian process models with each iteration $s$.  Eventually, given 
enough guesses and enough iterations, the stochastic model will break 
out of a local optimum to perhaps find a better global optimum elsewhere.
For the \emph{in situ} use case, though, eventual solution is not fast
enough.  The inner product search spans regions that the stochastic weighting
labeled as unlikely.  This assist is what enables the \allflag method
to converge faster than the stochastic method on its own.  That said, 
if the \emph{in situ} use case is not relevant, the stochastic method on 
its own could be used without any algorithmic qualms.

\section{Conclusions \& Future Work}
\label{conclusion}

The WORG method provides a deployment schedule optimizer that converges 
both closely enough and fast enough to be used inside
of a nuclear fuel cycle simulator. The algorithm can consistently obtain
tolerances of half-a-percent to a percent (1 GWe distances for over 200 GWe
deployable) for the once-through fuel cycle featured here within only five to 
ten simulations. Such optimization problems are made
more challenging due to the integral nature of facility deployment and
that any demand curve may be requested.

WORG works by setting up a Gaussian process to model the production
as a function of time and the deployment schedule. This model may then
be evaluated orders of magnitude faster than running a full simulation, enabling
the search over many potential deployment schedules. The quality of these
possible schedules is evaluated based on the dynamic time warping distance
to the demand curve. The lowest distance curve is then evaluated in a
full fuel cycle simulation. The production curve that is computed by the 
simulator in turn goes on to update the Gaussian process model and the
cycle repeats until the limiting conditions are met.

However, choosing the deployment schedules to estimate with the Gaussian
process may be performed in a number of ways. A blind approach would
simply be to choose such schedules randomly from a univariate. However, 
the WORG method has more information available to it that helps drive 
down the number of loop iterations. The first method discussed remains 
stochastic but uses the inverse DTW distances of the GP model to 
weight the deployment options, falling back to a Poisson distribution as 
necessary. This second method minimizes the model distance for each point 
in time from start to end, iteratively building up a solution. Finally, 
another estimation strategy tries both previous options and chooses the 
best result, forcing the stochastic method two of every four iterations 
to avoid deterministic loops.  It is this last all-of-the-above method 
that is seen to converge the fastest and to the lowest distance in most 
cases.

It is important to note that the WORG algorithm is applicable to any 
demand curve type and fuel cycle facility type. It is not restricted to 
reactors and power.  Enrichment and separative work units, reprocessing
and separations capacity, and deep geologic repositories and their
space could be deployed via the WORG method for any applicable demand 
curve.  Reactors were chosen for study here as the representative keystone 
example.

The next major step for this work is to actually employ the WORG method in 
a fuel cycle simulator.  However, to the best knowledge of the 
author, no existing simulator is capable of spawning forks of itself 
during run time, rejoining the processes, and evaluating the results of the 
child simulations in the parent simulation. Concisely, while many simulators 
are 
`dynamic' in the fuel cycle sense, none are `dynamic' in the programming
language sense. This latter usage of the term is what is required to 
take advantage of any sophisticated \emph{in situ} deployment optimizer.
The Cyclus fuel cycle simulator looks most promising as a platform
for such work to be undertaken. However, many technical roadblocks 
on the software side remain, even for Cyclus.

Furthermore, adding \emph{in situ} capability also adds the additional 
degree of freedom of how often to run the deployment schedule optimizer.
Running WORG each and every
time step seems excessive \emph{a priori}. Is every year, five years,
or ten years sufficient? How does this degree of freedom balance with the
time horizon $T$ specificed in the optimizer? These questions remain unanswered, even 
in a heuristic sense, and thus the frequency of optimization will be a key 
parameter in a future \emph{in situ} study.


\bibliographystyle{ans}
\bibliography{refs}

\begin{thebibliography}{10}

\bibitem{DBLP:journals/corr/HuffGCFMOSSW15}
K.~D. HUFF et~al.,
\newblock ``Fundamental Concepts in the Cyclus Fuel Cycle Simulator
  Framework,''
\newblock {\em CoRR}, {\bf abs/1509.03604} (2015).

\bibitem{cyclus_v1_0}
R.~W. CARLSEN et~al.,
\newblock ``{Cyclus v1.0.0},''
\newblock (2014),
\newblock http://dx.doi.org/10.6084/m9.figshare.1041745.

\bibitem{muller}
M.~M\"ULLER,
\newblock ``Dynamic Time Warping,''
\newblock in {\em Information Retrieval for Music and Motion}, pages 69--84,
  Springer Berlin Heidelberg, 2007.

\bibitem{rasmussen2006gaussian}
C.~E. RASMUSSEN and C.~K. WILLIAMS,
\newblock {\em Gaussian processes for machine learning},
\newblock The MIT Press (2006).

\bibitem{kelton2000simulation}
W.~D. KELTON and A.~M. LAW,
\newblock {\em Simulation modeling and analysis},
\newblock McGraw Hill Boston (2000).

\bibitem{vanderbei2001linear}
R.~J. VANDERBEI,
\newblock ``Linear programming,''
\newblock {\em Foundations and Extensions. Second Edition-International Series
  in Operations Research and Management Science}, {\bf 37} (2001).

\bibitem{kennedy2010particle}
J.~KENNEDY,
\newblock ``Particle swarm optimization,''
\newblock in {\em Encyclopedia of Machine Learning}, pages 760--766, Springer,
  2010.

\bibitem{vaz2009pswarm}
A.~I.~F. VAZ and L.~N. VICENTE,
\newblock ``PSwarm: A hybrid solver for linearly constrained global
  derivative-free optimization,''
\newblock {\em Optimization Methods \& Software}, {\bf 24}, 669 (2009).

\bibitem{osborne2009gaussian}
M.~A. OSBORNE, R.~GARNETT, and S.~J. ROBERTS,
\newblock ``Gaussian processes for global optimization,''
\newblock {\em Proc. 3rd international conference on learning and intelligent
  optimization (LION3)}, pages 1--15, 2009.

\bibitem{simpson2001kriging}
T.~W. SIMPSON, T.~M. MAUERY, J.~J. KORTE, and F.~MISTREE,
\newblock ``Kriging models for global approximation in simulation-based
  multidisciplinary design optimization,''
\newblock {\em AIAA journal}, {\bf 39}, 2233 (2001).

\bibitem{paciorek2004nonstationary}
C.~PACIOREK and M.~SCHERVISH,
\newblock ``Nonstationary covariance functions for Gaussian process
  regression,''
\newblock {\em Advances in neural information processing systems}, {\bf 16},
  273 (2004).

\bibitem{scikit-learn}
F.~PEDREGOSA et~al.,
\newblock ``Scikit-learn: Machine Learning in {P}ython,''
\newblock {\em Journal of Machine Learning Research}, {\bf 12}, 2825 (2011).

\bibitem{hodlr}
S.~{Ambikasaran}, D.~{Foreman-Mackey}, L.~{Greengard}, D.~W. {Hogg}, and
  M.~{O'Neil},
\newblock ``{Fast Direct Methods for Gaussian Processes and the Analysis of
  NASA Kepler Mission Data},''
\newblock (2014).

\end{thebibliography}
\end{document}